

\documentstyle{amsppt}
\magnification=\magstep1
\NoRunningHeads

\vsize=7.4in


\def\bull{\vrule height .9ex width .8ex depth -.1ex}


\topmatter
\title
An elementary example of a Banach space not isomorphic
to its complex conjugate
\endtitle

\author
N.J. Kalton \\
University of Missouri-Columbia
\endauthor

\address
Department of Mathematics,
University of Missouri
Columbia, MO  65211, U.S.A.
\endaddress
\email mathnjk\@mizzou1.bitnet \endemail
\thanks Supported by NSF grant DMS-9201357
 \endthanks
\subjclass
46B03
\endsubjclass
\abstract

We give a simple and explicit example of a complex Banach
space which is not isomorphic to its complex conjugate, and
hence of two real-isomorphic spaces which are not
complex-isomorphic.

\endabstract

\endtopmatter

\document \baselineskip=14pt

If $X$ is a complex Banach space, its complex conjugate
$\bar X$ is defined to be the space $X$ equipped with the
alternative scalar multiplication $\alpha\otimes
x=\bar\alpha x.$ In \cite{1} Bourgain gave the first example
of a complex Banach space which is not isomorphic to its
complex conjugate, thereby also producing an example of two
Banach spaces which are isomorphic as real Banach spaces,
but not isomorphic as complex Banach spaces; see also the
work of Szarek \cite{5} and Mankiewicz \cite{4}.  The
approach of all these authors was probabilistic, and no
explicit example has been constructed.  The aim of this note
is to give an elementary explicit construction of a Banach
space $X$ so that $X$ and $\bar X$ are non-isomorphic.  Our
example has the further advantage of being, in a certain
sense, entirely natural.

The space $X$ we construct is a {\it twisted Hilbert space},
i.e.  $X$ has a closed subspace $E$ so that $E$ and $X/E$
are Hilbertian; this further implies $X$ is hereditarily
Hilbertian.

Let $\omega$ denote the space of all complex-valued
sequences, $x=(x(n))_{n=1}^{\infty}$.  We let $e_n$ be the
canonical basis vectors.  Suppose $f:[0,\infty)\to\bold C$
is a Lipschitz map, with $f(0)=0.$ We define a map
$\Omega:\ell_2\to \omega$ by
$$\Omega_f(x)(n)=x(n)f(\log\frac{\|x\|_2}{|x(n)|}).$$ Here
we interpret the right-hand side to be zero if $x(n)=0.$ We
then define $Z_2(f)$ to be the space of pairs
$(x,y)\in\ell_2\times\omega$ so that $$ \|(x,y)\|_f =\|x\|_2
+\|y-\Omega_f(x)\|_2<\infty.$$ It then follows that $Z_2(f)$
is a Banach space under a norm equivalent to the quasi-norm
$\|\,\|_f.$ Such spaces were first considered in \cite{2},
where only the real case was considered; however, the switch
to complex scalars, and complex-valued $f$ is routine.

If $s\in \ell_{\infty}$ with $\|s\|_{\infty}\le 1$ then we
have the estimate $\|\Omega_f(sx)-s\Omega_f(x)\|_2 \le
C_0\|x\|_2$ where $C_0$ depends on the Lipschitz constant of
$f$ and this leads to fact that there is a constant $C$ such
that $\|(sx,sy)\|_f \le C\|(x,y)\|_f.$ It follows that the
spaces $E_n$ spanned by $(e_n,0)$ and $(0,e_n)$ form a UFDD
for $Z_2(f).$ This UFDD has a certain symmetry, which will
be used frequently, for if $\pi:\bold N\to\bold N$ is a
permutation and $x_{\pi}(n)=x(\pi(n))$ then
$\|(x_{\pi},y_{\pi})\|_f=\|(x,y)\|_f.$

We will now specialize to the functions
$f_{\alpha}(t)=t^{1+i\alpha}$ for $-\infty<\alpha<\infty.$
We write $Z_2(\alpha)=Z_2(f_{\alpha})$ and $\Omega_{\alpha}$
in place of $\Omega_{f_{\alpha}}.$

The following observation is trivial:

\proclaim{Proposition 1}The complex conjugate of
$Z_2(\alpha)$ is isomorphic to $Z_2(-\alpha).$ \endproclaim

\proclaim{Theorem 2}Suppose $Z_2(\alpha)$ and $Z_2(\beta)$
are isomorphic.  Then $\alpha=\beta.$\endproclaim

\demo{Proof}We suppose that $\alpha\neq 0$ and that
$Z_2(\alpha)$ is isomorphic to $Z_2(\beta).$ Let
$a=1+i\alpha$ and $b=1+i\beta.$

We first observe the following inequalities for $t>s\ge 0:$
$$ \align |t^b-s^b| &\ge t-s\tag1\\ |t^b-s^b| &\le
|b|(t-s)\tag2\\ |t^b-s^b-bs^{b-1}(t-s)| &\le \frac12
|b||s|^{-1}(t-s)^2.\tag3 \endalign $$

For $w\in \ell_2$ define
$\Omega'_{\beta}(w)(n)=w(n)(\log|w(n)|^{-1})^b.$ Note that
$$ \|\Omega_{\beta}(w)-\Omega'_{\beta}(w)\|_2 \le
|b||\log\|w\|_2|\|w\|_2.$$ In particular if $\|w\|_2\le 1$
$$ \|\Omega_{\beta}(w)-\Omega'_{\beta}(w)\|_2 \le |b|.$$ If
$A$ is a finite subset of $\bold N$ we will let
$\xi_A=\sum_{n\in A}e_n.$

We will suppose the existence of an operator $T:Z_2(\alpha)
\to Z_2(\beta)$ such that $\|T\|<1$ and, $c>0$ so that for
every $n,$ $\|T(e_n,0)\|_{\beta},\|T(0,e_n)\|_{\beta}>c.$ We
will say that $T$ is admissible if it satisfies these
properties, for some $c>0.$ Clearly if there is an
admissible operator then, by blocking, we can find an
admissible operator $T$ and an increasing sequence of
integers $(p_n)_{n=0}^{\infty}$ so that for suitable
sequences $u,v,w,y\in\omega$ we have, setting
$B_n=\{p_{n-1}+1,\ldots,p_n\},$
$T(e_n,0)=(u\xi_{B_n},v\xi_{B_n})$ and
$T(0,e_n)=(w\xi_{B_n},y\xi_{B_n}).$ Here
$u\xi_{B_n}=\sum_{k\in B_n}u(k)e_k$ etc.  Henceforward we
consider only operators blocked in this way.

Our first objective is to show that we must have
$\lim_{n\to\infty}\|w\xi_{B_n}\|_2=0.$ Indeed if this is not
the case, then loss of generality we may suppose that $T$
satisfies $\|w\xi_{B_n}\|_2\ge \delta>0$ for all $n.$ For
any integer $N$, let $A=\{1,2\ldots,N\}.$ Then
$\|T(0,N^{-1/2}\xi_A)\|_{\beta}<1$ and hence $$
\|N^{-1/2}\sum_{k=1}^Ny\xi_{B_k}
-\Omega_{\beta}(N^{-1/2}\sum_{k=1}
^Nw\xi_{B_k})\|_2 \le 1.$$ Since
$\|y\xi_{B_k}-\Omega_{\beta}(w\xi_{B_k})\|_2 \le 1, $ we
have $$
\|\Omega_{\beta}(N^{-1/2}\sum_{k=1}^N
w\xi_{B_k})-N^{-1/2}\sum_{k=1}^N
\Omega_{\beta}(w\xi_{B_k})\|_2 \le 1.$$ From this we get,
since $\|w\xi_{B_k}\|_2\le 1,$ $$
\|\Omega'_{\beta}(N^{-1/2}\sum_{k=1}^N
w\xi_{B_k})-N^{-1/2}\sum_{k=1}^N
\Omega'_{\beta}(w\xi_{B_k})\|_2 \le 1+2|b|.$$ This
simplifies to $$ \frac1N\sum_{k=1}^N\sum_{j\in B_k}
|w_j|^2|(\log(\frac{\sqrt N}{|w(j)|})^b
-(\log(\frac1{|w(j)|})^b|^2 \le (1+2|b|)^2.$$ From this we
deduce, using (1), that $$ \frac1N\sum_{k=1}^N\sum_{j\in
B_k} |w_j|^2(\log{\sqrt N})^2\le (1+2|b|)^2$$ which in turn
implies that $$ \delta^2(\log{\sqrt N})^2 \le (1+2|b|)^2$$
for all $N,$ leading to a contradiction.  We may thus
suppose that $\lim_{n\to\infty}\|w\xi_{B_n}\|_2=0.$

It now follows, by passing to a subsequence and rearranging
that we can further suppose that $w=0.$ We then have $c\le
\|y\xi_{B_n}\|_2 \le 1$ for all $n,$ and some $c>0.$

We next show again by contradiction that we cannot have
$\inf\|u\xi_{B_n}\|_{\infty}=0.$ Indeed, in the contrary
case, we can assume that for any $M$ there exists an
admissible $T=T_M$ (with $w=0$) so that $\|u\|_{\infty}\le
e^{-M}.$

Under this assumption suppose that $n_1<N_1$ and $n_2<N_2$
are pairs of integers and let $\sigma_r=\frac12\log n_r$ and
$\tau_r=\frac12\log N_r$ for $r=1,2.$ Let $N$ be any integer
greater than $\max(N_1,N_2).$ Suppose $A$ is a subset of
$\{1,2,\ldots,N\}$ with $|A|=n_1.$ Then
$\|(\xi_A,\sigma_1^a\xi_A)\|_{\alpha}=n_1^{1/2}.$ It follows
that $$ \|n_1^{-1/2}\sigma_1^a\sum_{k\in
A}y\xi_{B_k}+n_1^{-1/2}\sum_{k\in
A}v\xi_{B_k}-\Omega'_{\beta}(n_1^{-1/2}(\sum_{k\in
A}u\xi_{B_k}))\|_2 \le 1+|b|.$$ For ease of notation let
$\rho_j=-\log|u(j)|$ when $|u(j)|\neq 0$ and let $\rho_j=M$
otherwise; note that $\rho_j\ge M.$ We can rearrange the
preceding equation as $$ \frac1{n_1}\sum_{k\in A}\sum_{j\in
B_k}|\sigma_1^ay(j)+v(j)-(\sigma_1+\rho_j)^bu(j)|^2 \le
(1+|b|)^2.$$ By averaging we obtain $$
\frac1{N}\sum_{k=1}^{N}\sum_{j\in B_k}
|\sigma_1^ay(j)+v(j)-(\sigma_1+\rho_j)^bu(j)|^2 \le
(1+|b|)^2.$$ Similarly we obtain $$
\frac1{N}\sum_{k=1}^{N}\sum_{j\in B_k}
|\tau_1^ay(j)+v(j)-(\tau_1+\rho_j)^bu(j)|^2 \le (1+|b|)^2.$$

At this point we use the estimate (3):  $$
|(s+\rho_j)^b-\rho_j^b-bs\rho_j^{b-1}|\le
\frac{|b|s^2}{2M}.$$ Thus $$
|(\tau_1+\rho_j)^b-(\sigma_1
+\rho_j)^b-b(\tau_1-\sigma_1)\rho_j^{b-1}|
\le \frac{|b|\tau_1^2}{M}.$$

We thus obtain by the triangle law that
$$\frac1{N}\sum_{k=1}^N\sum_{j\in
B_k}|\frac{\tau_1^a-\sigma_1^a}{\tau_1-\sigma_1}y(j)-
b\rho_j^{i\beta}u(j)|^2 \le \frac{(M^{-1}|b|\tau_1^2 +
4+4|b|)^2}{(\tau_1-\sigma_1)^2}$$ using the fact that
$\|u\xi_{B_k}\|_2\le 1$ for all $k.$

Using a similar calculation for $\sigma_2,\tau_2$ and using
the fact that $\frac1N\sum_{k=1}^N\sum_{j\in B_k}|y(j)|^2
\ge c^2$ we obtain that $$
\big|\frac{\tau_1^a-\sigma_1^a}{\tau_1-\sigma_1}-
\frac{\tau_2^a-\sigma_2^a}{\tau_2-\sigma_2}\big| \le
K_1+K_2$$ where $$ K_r=\frac1c\big(\frac{ M^{-1}|b|\tau_r^2
+ 4+4|b|}{\tau_r-\sigma_r}\big).$$

Now $M$ is arbitrary, and $\tau_r,\sigma_r$ are restricted
only to be of the form $\frac12\log m$ where $m\in\bold N.$
It thus follows that for any $\kappa>1$ such that
$\kappa^a\neq 1$ we must have that
$$\lim_{\sigma\to\infty}\frac{\kappa^a\sigma^{a}-\sigma^a}{\kappa
\sigma-\sigma}
=\lim_{\sigma\to\infty}
\sigma^{i\alpha}\frac{\kappa^a-1}{\kappa-1}$$
exists, which contradicts the fact that $\alpha=\Re a\neq
0.$

Our conclusion is then there exists an admissible $T$ so
that $\inf\|u\xi_{B_n}\|_{\infty}>0.$ Under these
circumstances we can apply the diagonalization procedure of
Proposition 1.c.8 of \cite{3} and a subsequence argument to
produce an operator $S:Z_2(\alpha)\to Z_2(\beta)$ with
$\|S\|<1$ and so that $S(e_n,0)=(\lambda e_n,\mu e_n)$ and
$S(0,e_n)=(0,\nu e_n)$ with $\lambda\neq 0.$ Let
$A=\{1,2\ldots,N\}.$ Then
$\|S(N^{-1/2}\xi_A,\sigma^aN^{-1/2}\xi_A)\|_2 <1$ where
$\sigma=\frac12\log N.$ Hence, $$
|\nu\sigma^a-\lambda(\sigma+\log|\lambda|)^b+\mu|\le 1.$$ As
this holds for all $N$ we must have $a=b$ or $\alpha=\beta,$
as required.  \bull\enddemo

\proclaim{Corollary 3}The space $Z_2(\alpha)$ is not
isomorphic to its complex conjugate when $\alpha\neq
0.$\endproclaim

 \enddocument